\documentclass[12pt]{amsart}
\usepackage{latexsym}
\usepackage{amsmath,amsfonts}
 \numberwithin{equation}{section}
\newtheorem{theo}{Theorem}[section]

\newtheorem{prop}[theo]{Proposition}

\newcommand{\rbar}{\overline{\mathbb R}}
\newcommand{\cl}{\operatorname{cl}}
\newcommand{\Rset}{\mathbb R}
\begin{document}

\title[Probabilistic Pompeiu-Hausdorf metric]{ Completeness with respect 
to the probabilistic Pompeiu-Hausdorff metric}

\author{  Stefan Cobza\c s}
\address{\it Babe\c s-Bolyai University, Faculty of Mathematics
and Computer Science, Ro-3400 Cluj-Napoca, Romania,
 E-mail: scobzas@math.ubbcluj.ro}

\begin{abstract}
The aim of the present paper is to prove that the family of all  closed nonempty subsets
of a complete probabilistic metric space $L$ is complete with respect to the
probabilistic Pompeiu-Hausdorff  metric $H$. The same is true for the families of all
closed bounded, respectively compact, nonempty  subsets of $L$. If $L$ is a complete
random normed space in the sense of \v Serstnev, then the family of all nonempty closed
convex subsets of $L$ is also complete with respect to $H$.

  The probabilistic Pompeiu-Hausdorff metric  was defined and
studied by R.J. Egbert, Pacific J. Math. {\bf 24} (1968), 437-455, in the case of Menger
probabilistic metric spaces, and by    R.M. Tardiff, Pacific J. Math. {\bf 65} (1976),
233-251, in general probabilistic metric spaces. The completeness with respect to
probabilistic Pompeiu-Hausdorff metric of the space of all closed bounded nonempty
subsets of some  Menger  probabilistic metric spaces was proved by J. Kolumb\' an and A.
So\' os, Studia Univ. Babes-Bolyai, Mathematica, {\bf 43} (1998), no. 2, 39-48, and
{\bf 46} (2001), no. 3, 49-66.

AMS 2000 MSC:  46S50, 54E70

Key words: probabilistic metric, random normed spaces, Pompeiu-Hausdorf metric
\end{abstract}
\maketitle
\section{Introduction}

The study of probabilistic metric spaces (PM spaces for short) was initiated by K.
Menger \cite{men42} and A. Wald \cite{wa43}, in connection with some measurements
problems in physics. The positive number expressing the distance between two points
$p,q$ of a metric space is replaced by a distribution function (in the sense of
probability theory) $F_{p,q}:\rbar \to [0,1]$, whose value $F_{p,q}(x)$ at the point
$x\in \Rset$ can be interpreted  as the probability that the distance between $p$ and
$q$ be less than $x$. Since then the subject developed in various directions, an
important one being that of fixed points in PM spaces. Important contributions to the
subject have been done by A.N. \v Serstnev and  the Kazan school of probability theory,
see \cite{sers62,sers63e,sers64,sers64b} and the bibliography in \cite{SS}.
\par
A clear and thorough presentation  of the results up to 1983 is given in the book by B.
Schweizer and  A. Sklar \cite{SS}. Beside this book, at the present there are several
others dealing with various aspects of analysis in probabilistic metric spaces and in
probabilistic normed spaces -- V. Istr\u a\c tescu \cite{Is}, I. Istr\u a\c tescu and
Gh. Constantin \cite{IC,ICo}, V. Radu \cite{Rad}, S.-S. Chang and Y. J. Cho \cite{Cho},
O. Had\v zi\' c \cite{H}, O. Had\v zi\' c and E. Pap \cite{HP}. In the present paper we
shall follow the treatise \cite{SS}.
\par %
The probabilistic Pompeiu-Hausdorff metric  on the family of nonempty closed subsets of
a PM space  was defined by Egbert \cite{egb68} in the case of Menger PM spaces, and by
Tardiff \cite{tar76} in general PM spaces (see also \cite[\S12.9]{SS}), by analogy with
the classical case. Sempi \cite{semp92} used   the probabilistic Pompeiu-Hausdorff
metric to prove the existence of a completion of a PM space. Some results have been
obtained also by Beg and Ali \cite{beg96}.  %
\par %
As it is well known, the family of nonempty closed bounded subsets of a complete metric 
space is complete with respect to the Pompeiu-Hausdorff distance (see, e.g., 
\cite[Chapter 1]{Hand}). The aim of the present paper is to prove the probabilistic 
analogue of this result for the family of all nonempty closed  subsets of a 
probabilistic metric space. We shall prove that the families of all nonempty closed 
bounded, respectively  compact, subsets of a complete probabilistic metric space $L$ are 
 also complete with respect to the probabilistic Pompeiu-Hausdorff metric. 
If $L$ is a complete random normed space in the sense of \v Serstnev, then   the family 
of all nonempty closed convex subsets of $L$ is complete with respect to the 
Pompeiu-Hausdorf metric too. In the case of  Menger PM spaces 
$(L,\rho,\operatorname{Min}),\, $ and  $(L,\rho, W),\,$ with $t$-norms  
$\operatorname{Min}(s,t) = \min\{s,t\},\; s,t \in [0,1],\,$ respectively $W(s,t) = 
\max\{s+t-1,0\},\,$ the completeness of the space of all closed bounded nonempty subsets 
of $L$ with respect to the probabilistic Pompeiu-Hausdorff metric was proved by Kolumb\' 
an and So\' os in \cite{kol-sos98} and \cite{kol-sos01}.   In the case of a Menger PM 
space $(L,\rho,\operatorname{Min}),\, $ they proved also in \cite{kol-sos98} the 
completeness 
 of the family of all compact nonempty subsets of $L$.
 These completeness results were applied in \cite{kol-sos98,kol-sos01,kol-sos02}
to prove the existence of invariant sets for finite families of contractions in PM
spaces of random variables ($E$-spaces in the sense of Sherwood \cite{she69}, or
\cite[Ch. 9, Sect. 1]{SS}). %
\par %
 As in Aubin's book \cite{Aub}, I have adopted the term Pompeiu-Hausdorf metric.
 For a short comment on this fact, as well as on the similar case of
  the Painlev\' e-Kuratowski convergence for sequences of sets, see \cite[page xiv]{Aub}.
  \section{Preliminary notions}
Denote by $\Delta$ the set of {\it distribution functions},
meaning nondecreasing, left continuous functions $F:\rbar\to
[0,1]$ with $F(-\infty )=0$ and $F(\infty) =1.$  Let $D$ be the
subclass of
$\Delta$ formed by all functions $F\in \Delta$ such that %
\[ %
\lim_{x\to -\infty}F(x)= 0\quad\mbox{ and}\quad \lim_{x\to\infty}F(x) =1.
\] %
\par %
The {\it weak convergence } of a sequence $(F_n)$ in $\Delta$ to
$F\in\Delta$, denoted by $F_n\xrightarrow{w}F$, means that the equality%
\begin{equation}\label{eq2.1}
\lim_{n\to\infty}F_n(x) = F(x)
\end{equation}
holds for every continuity point $x$ of $F$. Since $F$ is
 non-decreasing the set of its discontinuity points is at most
 countable, so that the set of continuity points of $F$ is dense
 in $\Rset$.   In order that  $F_n\xrightarrow{w}F$
 it is sufficient that the relation \eqref{eq2.1} holds for every $x$ in an
 arbitrary dense subset of $\Rset$. An important result
 concerning weak convergence of distribution functions
 is Helly's First Theorem: every sequence in $\Delta$ contains a
 weakly convergent subsequence  (see Lo\` eve \cite[Sect. 11.2]{Lo}).
\par %
 The topology of weak convergence in $\Delta$ is metrizable. The
 first who realized this was P. L\' evy (see the Appendix to Fr\'
 echet's book \cite{Fr}), and for this reason  the metrics
generating  the weak convergence  in $\Delta$ are called L\' evy
 metrics. Since the original L\' evy metric characterizes the weak
 convergence only in $D$, Sibley \cite{sib} proposed a
 modification of L\' evy metric that generates the weak convergence
 in $\Delta.$ We shall work with a further modification proposed
 by Schweizer and Sklar \cite{SS} and denoted by $d_L$. The
distance $d_L(F,G)$  between two functions $F,G \in \Delta$ is
defined as the infimum of all numbers $h>0$ such that the
 inequalities %
\[ %
 F(x-h)-h\leq G(x) \leq F(x+h)+h
 \]%
 and %
 \[ %
G(x-h)-h\leq F(x)\leq G(x+h)+h
\] %
hold  for every $x\in (-h^{-1};h^{-1})$. One shows that $d_L$ is a metric on $\Delta$
and, for any sequence $(F_n)$
in $\Delta$ and $F\in \Delta$, we have %
\[ %
 F_n\xrightarrow{w}F\;\iff\; d_L(F_n,F)\to 0.
 \] %
By Helly's First Theorem the space $(\Delta,d_L)$ is compact,
hence complete (see \cite[\S 4.2]{SS}).
\par %
The sets of {\it distance functions } are:%
\[ %
 \Delta^+ = \{F\in \Delta : F(0) =0\} \quad\mbox{and}\quad D^+=D\cap \Delta^+.
 \] %
\par %
It follows that for $F\in \Delta^+$ we have  $F(x) = 0,\; \forall
x\leq 0 $. The set $\Delta ^+$ is closed in the metric space
$\Delta$, hence compact and complete too.
\par %
Two important distance functions are %
\[
\begin{aligned}
 \epsilon_0(x) &= 0\quad \mbox{for}\quad x\leq 0\qquad\qquad \qquad\mbox{and}\qquad \qquad
 \qquad
 &\epsilon_\infty(x) &= 0\quad \mbox{for}\quad x<\infty \\ \notag
               & = 1\quad\mbox{for}\quad x>0     &\;  & = 1\quad\mbox{for}\quad x=\infty
 \end{aligned}
 \]

\par %
The order in $\Delta^+$ is  defined  as the punctual order: for
$F,G\in \Delta^+$ we put %
\[ %
 F\leq G\quad\iff \quad \forall x >0\;\; F(x)\leq G(x).
 \] %
\par %
It follows that $\epsilon_0$ is the maximal element of $\Delta^+$
and of $D^+$ as well, and $\epsilon_\infty$ is the minimal element
of $\Delta^+$.
\par %
In the following we shall define some functions, say $F$,  on $\Rset$ and consider them
automatically extended to $\rbar$ by $F(-\infty)=0$ and $F(\infty) = 1.$
\par %
 If $\{F_i : i\in I\}$ is a family of functions in $\Delta^+$  then the function
 $F:\rbar \to [0,1]$ defined by %
 \[ %
 F(x) = \sup \{F_i(x) : i\in I\},\; x\in \Rset,
 \] %
 is the supremum of the family $\{F_i\}$ in the ordered set $(\Delta^+,\leq)\;$ -- \;
 $F=\sup_{i\in I}F_i.$
\par %
 To define the infimum of  the family $\{F_i\}$ put %
\begin{equation}\label{eq2.2}
\Gamma (x) = \inf\{F_i(x) : i\in I\},\; x\in \Rset. %
\end{equation}
\par %
 Since the function $\Gamma$ is nondecreasing, but not necessarily  left continuous on
 $\Rset$, we have to regularize it by taking the left limit %
 \begin{equation}\label{eq2.3}
 G(x) = \ell^-\Gamma(x) := \lim_{x'\nearrow x}\Gamma (x') =
 \sup _{x'<x}\Gamma (x'), \; x\in \Rset. %
 \end{equation} %
\par %
Then $G(x) \leq \Gamma(x),\; \forall x\in R,$ the function $G$ belongs to $\Delta^+$ and
$G=\inf_{i\in I}F_i$ -- the infimum of the family $\{F_i\}$ in the ordered set
$(\Delta^+,\leq)$.
\par %
A {\it triangle function} is a binary operation $\tau$ on $\Delta^+, \; \tau :
\Delta^+\times \Delta^+ \to \Delta^+,\;$ that is commutative, associative,
non-decreasing in each place ($\tau(F_1,G_1)\leq\tau(F_2,G_2),\;$ if $F_1\leq F_2$ and
$G_1\leq G_2$), and has $\epsilon_0$ as identity:  $\tau(F,\epsilon_0) = F,\; F\in
\Delta^+.$ The triangle function $\tau$ is called continuous if it is continuous with
respect to the $d_L$-topology of $\Delta^+$. It follows that $\tau$ is, in fact,
uniformly continuous, since the metric space $(\Delta^+,d_L)$ is compact.
\section{Probabilistic metric spaces} %
A {\it probabilistic metric space } (PM space) is a triple $(L,\rho,\tau)$, where $L$ is
a set, $\rho$ is a mapping from $L\times L$ to $\Delta^+$, and  $\tau$ is a continuous
triangle function.  The value of $\rho$ at $(p,q)\in L\times L$ is denoted by
$\,F_{pq},\;$ i.e., $ \; \rho(p,q) = F_{pq}.$
\par %
One supposes that the following conditions are satisfied for all $p,q,r \in L$:
\par %
(PM1) \quad $F_{pp} = \epsilon_0,$
\par %
(PM2) \quad $ F_{pq}=\epsilon_0\; \Rightarrow \; p=q,$
\par %
(PM3) \quad $F_{pq}=F_{qp},$
\par %
(PM4) \quad $ F_{pr}\geq \tau(F_{pq},F_{qr}).$
\par %
The mapping $\rho$ is called the {\it probabilistic metric} on $L$
and the condition (PM4) is the probabilistic analogue of the
triangle inequality.
\par %
The {\it strong topology} on a PM space is defined by the neighborhood system: %
\begin{equation}\label{eq3.1} %
 U_t(p) = \{q\in L : F_{pq}(t)> 1-t\}, \quad t>0. %
 \end{equation} %
\par %
 Putting %
 \begin{equation}\label{eq3.2}%
 \bar{U}_t(p) = \{q\in L : F_p(t) \geq 1-t\}
 \end{equation}
 we have $U_t(p) \subset \bar{U}_t(p)\;$ and $\; \bar{U}_{t'}(p)\subset U_t(p)\;$
 for $t'<t,$ showing that   the family \eqref{eq3.2} of subsets  of $L$
forms also a neighborhood base for the strong  topology of $L$.
\par %
  Observe that $U_t(p) = L,\; $ for $t > 1$, and $\bar{U}_t(p) = L, $ for $t\geq 1,$
so that we can restrict to $t\in (0,1)$ when working with strong neighborhoods. In fact,
we can suppose that $t$ is as small as we need.
\par %
 The strong topology on a PM space $(L,\rho,\tau)$ is derived from the uniformity
 $\mathcal U$ generated by the vicinities:  %
 \begin{equation}\label{eq3.2a} %
 U_t = \{(p,q)\in L\times L : F_{pq}(t)>1-t\}, \quad t>0.
 \end{equation} %
 The strong topology is metrizable since $\{U_{1/n} : n\in \mathbb N\}$ is a countable
 base for the uniformity $\mathcal U$.  The probabilistic metric
  $\rho$ is uniformly continuous mapping from $L\times L$ with the product topology to
  $(\Delta^+,d_L)$,  meaning that %
   \begin{equation}\label{eq3.3}%
  p_n\to p\;\;\mbox{and}\;\; q_n \to q \;\; \mbox{in} \;\;L\;\; \Rightarrow
  \;\;
   F_{p_{n}q_{n}}\xrightarrow{w} F_{pq}.
\end{equation}
\par %
  The convergence of  a sequence $(p_n)$ in $L$ to $p\in L$ is characterized by %
\[ %
\begin{aligned} %
 p_n\to p \; \iff &\; \forall t>0 \; \exists n_0 \;\forall n\geq n_0\quad p_n\in U_t(p)\\
             \iff & \; F_{p_np}\xrightarrow{w} \epsilon_0 \\
             \iff & \; d_L(F_{p_np},\epsilon_0) \to 0.
\end{aligned}
\]
\par %
A sequence $(p_n)$ in $L$ is called a {\it Cauchy sequence}, or { \it fundamental}, if %
\[ %
 F_{p_np_m}\xrightarrow{w} \epsilon_0 \quad\mbox{for} \quad n,m\to\infty,
 \] %
or, equivalently, %
\[ %
\forall t>0 \; \exists n_0\; \mbox{such that}\; \forall n,m \geq
n_0 \; \; (p_n,p_m)\in U_t \; (\iff F_{p_{n}p_{m}}(t) > 1-t).
\] %
\par %
A convergent sequence in $L$ is a Cauchy sequence, and the PM space $L$ is called {\it
complete} (with respect to the strong topology) if every Cauchy sequence is convergent.
\par %
For these and other questions concerning the strong topology of a PM space, %
see \cite[Chapter 12]{SS}.
\par %
Throughout this paper all the topological notions concerning a PM
space will be considered with respect to the strong topology.
\section{The probabilistic Pompeiu-Hausdorff metric} %
For  a metric space $(X,d)$, two nonempty bounded subsets $A,B$ of
$X$ and a point $p\in
X$, one introduces the following  notations and notions : %
\[ %
\begin{aligned} %
d(p,B) &= \inf\{d(p,q) : q\in B\}\; - \; \mbox{ the distance from}\; p \;\mbox{to}\; B,\\ %
 h^*(A,B) &= \sup\{d(p,B) : p\in A\}\;-\; \mbox{the excess of}\; A \; \mbox{ over} \; B,
 \end{aligned}
\] %
 and let %
 \[ %
 h(A,B) = \max \{h^*(A,B),h^*(B,A)\}\;
 \] %
 be the Pompeiu-Hausdorff distance between the sets $ A,B$.
\par %
 Denoting by $P_{fb}(X)$ the family of all nonempty closed bounded subset of $X$ it follows that
 $h$ is a metric on $P_{fb}(X)$, and the metric space $(P_{fb}(X),h)$  is complete if
 $(X,d)$ is complete (see, e.g., \cite[Chapter 1]{Hand}).
\par %
 In the case of a PM space $(L,\rho,\tau)$ the definitions are similar but, taking into
 account the fact that the probabilistic triangle inequality (PM4)  is written in
 reversed form with respect to the usual triangle inequality,  sup and inf will change their
 places.
\par %
 For   two nonempty subsets  $A,B$ of  $ L$ and $p\in L$ denote by %
 \begin{equation}\label{eq4.1} %
 F_{pB} = \sup\{ F_{pq} : q\in B\} \; \iff \;
  F_{pB}(x) = \sup\{ F_{pq}(x) : q\in B\},\; x\in \Rset,
\end{equation} %
the {\it probabilistic distance } from $p$ to $B$, and let %
  \begin{equation}\label{eq4.2}  %
F^*_{AB} = \inf\{F_{pB} : p\in A\}.
\end{equation}
\par %
Taking into account the formulae \eqref{eq2.2} and \eqref{eq2.3}, it follows %
\[ %
 F^*_{AB} = \ell^-\Gamma^*_{AB},
 \] %
where %
\[ %
 \Gamma^*_{AB}(x) = \inf \{F_{pB}(x) : p\in A\}, \; x\in \Rset.
 \] %
\par %
The {\it probabilistic Pompeiu-Hausdorff distance}  between  the sets $A,B$ is defined
by  \newline %
 $\;H(A,B)\; =   F_{AB},\;$  where %
\begin{equation}\label{eq4.3} %
F_{AB}(x) = \min\{ F^*_{AB}(x), F^*_{BA}(x)\},\;\; x\in \Rset.
\end{equation} %
The probabilistic Pompeiu-Hausdorff metric was defined and studied  by Egbert
\cite{egb68} in the case of Menger PM spaces and by Tardiff \cite{tar76} in general PM
spaces (see also \cite[\S 12.9]{SS}). The  mapping $H(A,B) = F_{AB}$ satisfies the
following properties, where
$\cl$ denotes the closure with respect to the strong topology: %
\begin{prop}\label{p4.1} {\em (\cite[Th. 12.9.2]{SS})}
\par %
1. \quad $F_{\{p\}\{q\}} = F_{pq}$ \quad for \; $p,q\in L;$
\par %
2. \quad For nonempty \;$A,B \subset L,\; F_{AB} = F_{BA}, \;\; F_{AB} =
F_{\cl(A)\cl(B)}$,\; and \; $F_{AB} = \epsilon_0 \;$ if and only if \; $\cl(A) =
\cl(B).$
\end{prop} %
In order that $H$ satisfy the probabilistic triangle inequality
(PM4), we have to impose a supplementary condition on the triangle
function $\tau$.
The triangle function $\tau$ is called {\it sup-continuous } if %
\begin{equation}\label{eq4.4} %
\tau(\sup_{i\in I} F_i,G) = \sup_{i\in I}\tau(F_i,G)
\end{equation}
for any family $\{ F_i : i\in I\} \subset \Delta^+\;$ of  distance functions and any
$G\in \Delta^+$.
\par %
Denote by  $P_{f}(L)$   the family of all nonempty closed subsets of a PM space
$(L,\rho,\tau)$.
\begin{theo}\label{th4.2}{\em (\cite[Th. 12.9.5]{SS})}  If the triangle function $\tau$
is sup-continuous then the mapping  $ H(A,B) = F_{AB},\,$ where $\,F_{AB}\,$ is defined
by  \eqref{eq4.3}, is a probabilistic metric on $P_f(L).$
\end{theo} %
In the following proposition we collect some properties which will
be used in the proof of the completeness of $P_f(L)$ with respect
to the probabilistic  Pompeiu-Hausdorff metric. %
\begin{prop}\label{p4.3}
Let $\;(L,\rho,\tau)\;$ be a PM space with sup-continuous triangle function $\tau$, and
let $A,B \in P_f(L)\;$ and $p\in L.$ Then
\par %
1. \quad $F_{pB} \geq \tau\left(F_{pA},F^*_{AB}\right)$;\\
and %
\par %
2. \quad $F_{pB} \geq \Gamma^*_{AB} \geq F^*_{AB} \geq F_{AB}.$ %
\par %
3. If $F_{AB}(s) > 1-s$ for some $s,\; 0<s<1,$ then  %
\begin{equation}\label{eq4.12}
\forall p\in A\; \exists q\in B \; \mbox{such that}\; F_{pq}(s)>1-s,
\end{equation}
and %
\begin{equation}\label{eq4.13} %
\forall q\in B\; \exists p\in A \; \mbox{such that}\; F_{pq}(s)>1-s.
\end{equation}
\end{prop}
\begin{proof}
 For $x\in \Rset$ we have  %
 $$
 \forall a\in A \;\forall b\in B\quad
 F_{pB}(x) \geq F_{pb}(x) \geq \tau(F_{pa}F_{ab})(x).$$
 Taking the supremum with respect to $b\in B$  and taking in account that $\tau $ is sup
 continuous and monotonic in each place, we get %
 $$
\forall a\in A \;\;  F_{pB}(x) \geq \tau(F_{pa},F_{aB})(x) \geq
\tau(F_{pA},F^*_{AB})(x). $$ %
Taking now the supremum with respect to $a\in A$ one obtains the
inequality 1.
\par %
The inequalities  2 are immediate from definitions.

To prove 3, observe that %
$$%
F_{AB}(s) > 1-s \; \iff \; F^*_{AB}(s)> 1-s \;\mbox{and}\; F^*_{BA}(s) > 1-s. $$
 It follows %
$$%
\inf \{ F_{p'B}(s): p' \in A\} = \Gamma^*(s) \geq F^*_{AB}(s) > 1-s, $$%
so that %
$$%
\sup \{F_{pq}(s) : q\in B\} > 1-s, $$%
implying \eqref{eq4.12}.

The inequality \eqref{eq4.13} can be proved similarly.
\end{proof} %

 The completeness result will be obtained under a further restriction imposed to
 $\tau$. We say that the triangle function $\tau$ satisfies the condition (W) if
$$%
 \mbox{(W)}\quad \quad F(x) > \alpha \;\mbox{and} \; G(x) > \beta \;
 \Rightarrow \; \tau(F,G)(x) > \max\{\alpha + \beta -1,0\},
$$ %
 for all $x>0,\;$ where  $F,G \in \Delta^+,\;$ and $\; \alpha, \beta \in \Rset.$
\par %
{\bf Remark.\;} Considering the $t$-norm %
$$ W(x,y) = \max \{x+y-1,0\}, \quad (x,y) \in [0,1]^2 , $$ %
(see \cite[p. 5]{SS}) and the associated triangle function $\,\bf W$,
 defined for $F,G \in \Delta^+$ by  %
$$ {\bf W}(F,G)(x) = W(F(x),G(x)), \; x\in \Rset,$$ %
(see  \cite[p. 97]{SS}),  the condition  (W) essentially means that $\tau \geq \bf W.$

\par %
Now we are ready to state and prove the completeness result.
\begin{theo}\label{th4.4}
Let $(L,\rho,\tau)$ be a PM space with sup-continuous triangle function $\tau$
satisfying the condition (W).
\par %
If the PM space $L$ is complete then the space $P_f(L)$ is complete with respect to the
probabilistic Pompeiu-Hausdorff metric.
\end{theo}
\begin{proof}
Let $(A_n)$ be a sequence in $P_f(L)$ that is fundamental with respect to the
probabilistic  Pompeiu-Hausdorff metric $H$.
\par %
Put %
\[%
 A = \bigcap_{n\geq 1}\cl\left(\bigcup_{m\geq n}A_m\right),
 \] %
and show that $A\in P_f(L)$  (meaning that $A\subset L$ is nonempty closed) and that
the sequence $(A_n)$ converges to $A$ with respect to the probabilistic
Pompeiu-Hausdorff metric $H$.
\par %
Observe that %
\begin{equation}\label{eq4.5} %
p\in A \quad \iff \quad \exists \;\; n_1<n_2<... \;\; \exists p_k\in A_{n_k}: \quad
p_k\to p. \end{equation} %
\par %
For $0< t < 1/2$  fixed, choose $n_0\in \mathbb N$ such that %
\[ %
 \forall n,m \geq n_0 \quad F_{A_nA_m}(t) > 1-t.
 \] %
\par %
For $m\geq n_0$ fixed, put $n_1 := m$ and pick an element $p_1\in A_{n_1}$.
\par %
Let now $n_2 > n_1$ be such that  %
\[ %
 \forall n,n' \geq n_2 \quad  F_{A_nA_{n'}}(\frac{t}{2}) > 1-\frac{t}{2}.
 \] %
The inequalities %
 \[ %
  F^*_{A_{n_1}A_{n_2}}(t) \geq F_{A_{n_1}A_{n_2}}(t) > 1-t
  \] %
 and the fact that $p_1$ belongs to $A_{n_1}$  imply $F_{p_1A_{n_2}}(t) > 1-t, $ so that
 there exists $p_2 \in A_{n_2}$ such that %
 \[%
 F_{p_1p_2}(t) > 1-t.
 \]%
 Take now $n_3 > n_2 $ such that %
\[%
 \forall n,n' \geq n_3 \quad F_{A_nA_{n'}}(\frac{t}{2^2}) > 1- \frac{t}{2^2}.
 \]%
 Reasoning like above, we can find an element $p_3\in A_{n_3}$ such that %
 \[%
  F_{p_2p_3}(\frac{t}{2})> 1-\frac{t}{2}.
  \]%
 Continuing in this way, we obtain a strictly increasing sequence of indices
 $n_1 < n_2 < ...$ and the elements $p_k \in A_{n_k},\; k\in \mathbb N, \;$ such that %
\begin{equation}\label{eq4.6} %
 F_{p_kp_{k+1}}(\frac{t}{2^{k-1}}) > 1 - \frac{t}{2^{k-1}}, %
 \end{equation}
 for all $k\in \mathbb N.$\\

 {\it Claim I.} \quad $ \forall i\in \mathbb N \; \forall k\in \mathbb
 N\quad  F_{p_kp_{k+i}}(\frac{t}{2^{k-1}}) >
 1- (\frac{1}{2^{k-1}}+\frac{1}{2^k}+...+\frac{1}{2^{k+i-2}})t.
 $ \\

 We proceed by induction on $i$.  For $i=1$ the assertion is true by the choice of the
 elements $p_k$ (see \eqref{eq4.6}).
\par %
 Suppose that the assertion is true for $i$ and prove it for $i+1$. Appealing to
 condition (W) we have %
 \[%
  F_{p_kp_{k+i+1}}(\frac{t}{2^{k-1}}) \geq \tau
 \left(F_{p_kp_{k+1}},F_{p_{k+1}p_{k+i+1}}\right)(\frac{t}{2^{k-1}}) >
  1-(\frac{1}{2^{k-1}}+\frac{1}{2^{k}}+...+ \frac{1}{2^{k+i-1}})t,
  \]%
  since $F_{p_kp_{k+1}}(\frac{t}{2^{k-1}}) > 1 -\frac{1}{2^{k-1}}$ and, by the induction
  hypothesis, %
  \[%
   F_{p_{k+1}p_{k+i+1}}(\frac{t}{2^{k-1}})\geq F_{p_{k+1}p_{k+i+1}}(\frac{t}{2^{k}})>
  1-(\frac{1}{2^{k}}+...+ \frac{1}{2^{k+i-1}})t.
  \]%
\vspace{1em}\\ %
 {\it Claim II. The sequence $(p_k)$ is
fundamental  in the PM space $L$.  }\\

For  $0< s < 1$ choose $k_0\in \mathbb N$ such that $2^{-k_0+1} < s$. Then for any
$k\geq k_0$ and arbitrary $i\in \mathbb N$ we have %
\[%
 F_{p_kp_{k+i}}(s) \geq F_{p_kp_{k+i}}(\frac{t}{2^{k-1}}) > 1-
(\frac{t}{2^{k-1}}+...+\frac{t}{2^{k+i-1}}) > 1-\frac{t}{2^{k}}
>1-s.
\]%
Since the PM space $L$ is complete, there exists $p\in L$ such
that $p_k\to p$ in the strong topology of $L$. The choice of the
elements $p_k$ and \eqref{eq4.5} yield $p\in A$. Since the set $A$
is obviously closed it follows $A\in P_f(L)$.
\par %
By Claim I we have %
\[%
 F_{p_1p_{k}}(t) > 1- (1+\frac{1}{2}+...+\frac{1}{2^{k-2}})t > 1-2t.
 \]%
Let now $t', \; t<t'<2t,\;$ be a continuity point of the
distribution
 function $F_{p_1p}$. The continuity of the distance function (see \eqref{eq3.3})
 and the inequalities %
\[%
 F_{p_1p_k}(t')\geq F_{p_1p_k}(t) > 1-2t %
\] %
yield, for $\;k\to \infty,\;\; F_{p_1p}(t')\geq 1-2t,\,$  so that %
\[%
 F_{p_1A}(t') = \sup_{q\in A} F_{p_1q}(t') \geq 1-2t.
 \]%
As $p_1$ was arbitrarily chosen in $A_m$, it follows %
\[ %
 \Gamma^*_{A_mA}(t') = \inf\{ F_{p'A}(t') : p'\in A_m\} \geq 1-2t.
 \]%
But then %
\[ %
 F^*_{A_mA}(2t) = \sup_{t'} \Gamma^*_{A_mA}(t')\geq 1-2t,
\] %
  where the supremum is taken
over all continuity points $t'$ of the function $F_{p_1p} $ lying in the interval
$(t,2t)$. The fact that the set of these points is dense in the interval  $(t,2t)$
justifies the  equality sign in the first of the above relations.
\par %
Taking into account that $m\geq n_0$ was arbitrarily chosen too, we finally obtain  %
\begin{equation}\label{eq4.7} %
 \forall m\geq n_0 \quad F^*_{A_mA}(2t)\geq 1-2t, %
\end{equation}
\par %
Let now $p\in A$ and let $n_1<n_2<...$ and $p_k\in A_{n_k}$ be such that $p_k\to p$ in
the strong topology of the PM space $L$.
\par %
Choose $k_0 \in \mathbb N$ such that %
\[%
 \forall k\geq k_0 \quad F_{p_kp}(t) > 1-t.
 \]%
\par %
 Proposition \ref{p4.3}, the inequality $F_{pA_{n_{k_0}}}\geq F_{pp_{k_0}}$,  and condition
(W) give, for any $t',\; t< t' < 2t,$ %
\[%
F_{pA_m}(t') \geq F_{pA_m}(t) \geq
\tau\left(F_{pA_{n_{k_0}}},F^*_{A_{n_{k_0}}A_m}\right)(t) \geq
\tau\left(F_{pp_{k_0}},F^*_{A_{n_{k_0}}A_m}\right)(t) > 1-2t. %
\]%
\par %
Since $p\in A$ was arbitrarily chosen, it follows %
\[%
\forall t',\; t<t'<2t,\quad \Gamma^*_{AA_m}(t') \geq 1-2t,
\] %
so that %
\begin{equation}\label{eq4.8} %
 \forall m\geq n_0 \quad F^*_{AA_m}(2t) \geq 1-2t. %
\end{equation}
 The inequalities \eqref{eq4.7} and \eqref{eq4.8} yield %
\[ %
\forall m\geq n_0 \quad H(A_m,A)(2t) = F_{A_mA}(2t) \geq 1-2t,
\] %
i.e., the sequence $(A_m)$ converges to $A$ with respect to the probabilistic
Pom\-peiu-Hausdorff metric $H$. \par %
 The proof of the completeness is complete.  \end{proof}

 The {\it diameter} of a subset  $A$ of a PM space $(L,\rho,\tau)$ is defined by %
 $$ D_A(t) = \ell^-\Phi_A(t) $$ %
 where %
  $$ \Phi_A(t) = \inf \{F_{pp'}(t) : p,p'\in A\}.$$
  The set $A$ is called {\it bounded} if $D_A \in D^+$, i.e. $\sup\{D_A(t) : t> 0\} =1\,$
  (see\cite[pages 200-201]{SS}).
  This is equivalent to %
  \begin{equation}\label{eq4.8a}
  \sup\{\Phi_A(t) :  t > 0\} = 1.
  \end{equation}

  Now we shall show that the families  $\,P_{fb}(L)\,$ and $\, P_k(L)\,$
  of all closed bounded nonempty subsets of a PM space $L$, respectively
  of all nonempty compact subsets of
$L$,  are complete in $P_f(L)$ with respect to the Pompeiu-Hausdorff metric, provided 
the  PM space $L$ is complete. To prove the assertion concerning the compact sets, we 
shall use the characterization of compactness in uniform spaces in terms of total 
boundedness (see \cite[Ch. 6]{Kel}). Let 
  $(X,\mathcal U)$ be a uniform space. For $U\in \mathcal U$  and a subset $A$ of $X$ put
 $U(A) = \{x\in X : \exists y\in A \; \mbox{such that} \; (x,y)\in U\}$. It follows that
 $U(x) = U(\{x\})$ is a neighborhood of $x$ and $\{U(x) : U\in \mathcal U\}$ forms a
 neighborhood base at $x$. A subset $Y$ of $X$ is called {\it totally bounded} if for
 every $U\in \mathcal U$ there exists a finite subset $Z$ of $X$ such that $Y\subset
 U(Z).$ Then a subset of a uniform space $(X,\mathcal U)$ is compact if and only
 if it is complete and totally bounded (\cite[Ch. 6, Th. 32]{Kel}). If $L$ is  a PM space then, 
considering  $L$ as a uniform space with respect to the uniformity 
 generated by the vicinities \eqref{eq3.2a}, 
  denote by $P_{ftb}(L)$ the family of all nonempty, closed and totally bounded
  subsets of $L$.  %
\begin{theo}\label{th4.5}
 If $(L,\rho,\tau)$ is a PM space with sup-continuous triangle function $\tau$
 satisfying the condition (W), then the subspaces $\,P_{fb}(L\,)$ and
 $\,P_{ftb}(L)\,$ are  closed in $P_f(L)$.

 Consequently, if the PM space $L$ is complete then the subspaces $\,P_{fb}(L)\,$ and
  $\,P_k(L)\,$ are complete with  respect to the probabilistic Pompeiu-Hausdorff metric.
\end{theo}
\begin{proof}
Let $(A_n)$ be  a sequence of closed bounded nonempty sets converging to $A\in P_f(L)$
with respect to probabilistic Pompeiu-Hausdorff metric $H$. We have to show that $A$ is
bounded too, i.e. that %
\begin{equation}\label{eq4.8a'} %
\sup\{\Phi_A(t) : t>0\} = 1. %
\end{equation} %

Let $0<\epsilon <1/3\,$ and let $m\in \mathbb N$ be such that %
\begin{equation}\label{eq4.8b} %
\forall n\geq m \quad F_{AA_n}(\epsilon) > 1-\epsilon. %
\end{equation}

Since $\sup\{\Phi_{A_m}(t) : t> 0\} = 1$ there exists $t > 0$ such that $\Phi_{A_m}(t) >
1-\epsilon,\,$ so that%
\begin{equation}\label{eq4.8c}
\forall q,q' \in A_m \quad F_{qq'}(t) > 1-\epsilon.
\end{equation}
We can suppose also that $t\geq \epsilon$. By \eqref{eq4.8b} and \eqref{eq4.12}, for any
$p,p'\in A$ there exist $q,q'\in A_m$ such that %
\begin{equation}\label{eq4.8d}
F_{pq}(\epsilon) > 1-\epsilon \quad\mbox{and}\quad F_{p'q'}(\epsilon) > 1-\epsilon. %
\end{equation}
Since $t\geq \epsilon$ we have $F_{pq}(t)\geq F_{pq}(\epsilon) > 1-\epsilon$ and
$F_{p'q'}(t)\geq F_{p'q'}(\epsilon) > 1-\epsilon,$  so that, by \eqref{eq4.8c}
and condition (W), we have %
$$ %
F_{qp'}(t) \geq \tau (F_{qq'},F_{q'p'})(t) > 1-2e,$$ %
and %
$$%
F_{pp'}(t) \geq\tau(F_{pq},F_{qp'})(t) > 1-3\epsilon. $$

We have proved that for any $\epsilon,\; 0<\epsilon <1/3,\;$ there exists $t>0$ such
that  $F_{pp'}(t) > 1-\epsilon\,$ for all $p,p'\in A$. It follows $\Phi_A(t)\geq 1-
3\epsilon,\, $ so that \eqref{eq4.8a'} holds.

 Suppose now that $(A_n)$ is a sequence of nonempty compact subsets of $L$
converging with respect to the probabilistic Pompeiu-Hausdorff metric $H$ to a set $A\in
P_f(L).$ We shall show that $A$ is totally bounded with respect to the uniformity having
as vicinities the sets $U_t$ given by \eqref{eq3.2a}.

Let $0<\epsilon <1/2$  and let $\,n\in \mathbb N$ be such that $F_{AA_n}(\epsilon)>
1-\epsilon.$ By \eqref{eq4.12} it follows %
\begin{equation}\label{eq4.9} %
\forall p\in A\; \exists q\in A_n\; \mbox{such that}\; F_{pq}(\epsilon) > 1-\epsilon. %
\end{equation} %

Now, since the set $A_n$ is totally  bounded, there exits a finite set $Z\subset L$ such
that %
\begin{equation}\label{eq4.10}
\forall q\in A_n \; \exists z\in Z \; \mbox{such that} \; F_{qz}(\epsilon ) > 1-
\epsilon . %
\end{equation}

For an arbitrary $ p\in A$ choose first an element $q\in A_n$ according to \eqref{eq4.9}
and then, for this $q$ select $z\in Z$ according to \eqref{eq4.10}. Taking into account
the condition (W) we get %
$$%
F_{pq}(\epsilon ) \geq \tau (F_{pq},F_{qz})(\epsilon ) > \max \{ 1-2\epsilon ,0\} =
1-2\epsilon . $$%
 It follows  $A\subset U_{2\epsilon}(Z)$, i.e. the set $A$ is totally bounded.

 Now, if  the PM space $L$ is complete and $A$ is closed in $L$, it follows that $A$ is
 complete too, hence  compact, as complete and totally bounded.
 \end{proof}

 {\bf Remark.\;} As we have yet mentioned, in the case of Menger PM spaces
$(L,\rho,\operatorname{Min}),\, $ and  $(L,\rho, W),\,$ the
completeness of the space of all closed bounded subsets of $L$ was
proved by Kolumb\' an and So\' os in \cite{kol-sos98} and
\cite{kol-sos01}, respectively. Since $\operatorname{Min}\geq W,\,
$ both of these results are contained in the above completeness
result.
  The completeness of $P_k(L)$ in the case of a Menger PM space
$(L,\rho,\operatorname{Min})$ was proved in \cite{kol-sos98}.

For a subset $A$ of a PM space $(L,\rho,\tau)$ and $0<\epsilon \leq 1$ let %
$$ %
A_\epsilon = \{q\in L : \exists p\in A\; F_{pq}(\epsilon)>1-\epsilon\} = \bigcup
\{U_\epsilon(p) : p\in A\}. $$%

As in the case of ordinary metric spaces we have: %
\begin{prop}\label{p4.6} %
(i)\quad $\cl A = \bigcap_{\epsilon >0}A_\epsilon$

If $\tau$ satisfies (W) then

(ii)\quad $A\subset B_\epsilon \; \Rightarrow \; \cl A \subset B_{2\epsilon }.$
\end{prop}
\begin{proof}
 Let $q\in \cl A$ and $\epsilon > 0.$ Choosing $p\in U_\epsilon(q)\cap A$ it follows
 $$%
 q\in U_\epsilon (p) \subset A_\epsilon. $$%
 i.e. $\cl A\subset \cap_\epsilon A_\epsilon.$
 To prove the reverse inclusion we shall show that %
 $$ \cap_{n\geq 1} A_{1/n} \subset \cl A. $$%
 If $q\in \cap_{n\geq 1} A_{1/n}$ then %
 $$ \forall n\; \exists p_n\in A \;\mbox{such that}\;  F_{pp_n} > 1-\frac{1}{n},$$ %
 which implies that $(p_n)$ converges to $p$ in the strong topology of the PM space $L$,
 i.e. $\; p\in\cl A.$

 To prove (ii), let $p\in \cl A$.  It follows $U_\epsilon (p) \cap A \neq \emptyset,\,$
 so that $F_{pq}(\epsilon) > 1-\epsilon,\;$ for some $\, q\in A.$

 Since $A\subset B_\epsilon$ it follows  $F_{qr}(\epsilon)>1-\epsilon,\;$ for some
 $r\in B$. But then, taking into account the condition (W), we have for $0<\epsilon \leq 1/2$%
 $$%
 F_{pr}(\epsilon)\geq \tau(F_{pq},F_{qr})(\epsilon) > \max \{1-2\epsilon , 0\} =
 1-2\epsilon, $$%
showing that $p\in B_{2\epsilon}$. If $\epsilon >1/2$ then $B_{2\epsilon} = L.$
\end{proof}

In the following proposition we give two expressions for the probabilistic
Pompeiu-Hausdorff limit of a sequence of sets in $\,P_f(L),\,$  inspired by a well known
result for the usual Pompeiu-Hausdorff metric (see \cite[Proposition 1.3]{Hand}).
\begin{prop}\label{p4.7}
Let $(L,\rho,\tau)$ be a PM space with sup-continuous triangle function $\tau$
satisfying the condition (W). If $(A_n)$ is sequence in $P_f(L) $ converging to $A\in
P_f(L)$ with respect to the probabilistic Pompeiu-Hausdorff metric $H$ then
\begin{equation}\label{eq4.11}
A =\bigcap _{n\geq 1}\cl \left(\bigcup_{m\geq n}A_m\right) = %
\bigcap_{\epsilon > 0}\bigcup_{n\geq 1}\bigcap_{m\geq n}(A_m)_\epsilon .
\end{equation}\end{prop}
\begin{proof}
   Show first that %
\begin{equation}\label{eq4.14}
A\subset \bigcap _{n\geq 1}\cl \left(\bigcup_{m\geq n}A_m\right).%
 \end{equation} %
 Let $p\in A$ and let $n_1\in \mathbb N$ be such that %
 $$%
 \forall m\geq n_1\; F_{AA_m}(\frac{1}{2})> 1- \frac{1}{2}.$$%
 By \eqref{eq4.12}, %
 $$%
 \exists p_1\in A_{n_1} \; F_{pp_1}(\frac{1}{2})>1-\frac{1}{2}.$$%
 Continuing in this way we obtain a sequence $n_1 < n_2 <...$ of indices and the
 elements $p_k\in A_{n_k}$ such that %
 $$%
  F_{pp_k}(\frac{1}{2^k}) > 1- \frac{1}{2^k}.$$ %
 It follows $p_k\to p$, so that %
 $$ p \in \bigcap _{n\geq 1}\cl \left(\bigcup_{m\geq n}A_m\right).$$
  Let now $0< \epsilon <1/2\,$ and let $n_0\in \mathbb N$ be such that %
  $$%
  \forall m\geq n_0 \; F_{AA_m}(\epsilon)>1-\epsilon.$$%
By \eqref{eq4.13} it follows %
$$%
\forall m\geq n_0 \; \forall q\in A_m \; \exists p\in A \;\mbox{such that} \;
F_{pq}(\epsilon) > 1 -  \epsilon, $$%
so that %
$$%
\forall m\geq n_0 \; A_m \subset A_\epsilon, $$ %
or, equivalently, %
$$%
\bigcup_{m\geq n_0}A_m \subset A_\epsilon.$$%
But then %
$$%
\bigcap _{n\geq 1}\cl \left(\bigcup_{m\geq n}A_m\right) \subset \cl \left(\bigcup_{m\geq
n_0}A_m\right) \subset A_{2\epsilon}. $$%
Since $0< \epsilon < 1/2$ is arbitrary we have %
$$%
\bigcap _{n\geq 1}\cl \left(\bigcup_{m\geq n}A_m\right) \subset \bigcap_{0<\epsilon <
1/2} A_{2\epsilon} = \cl A = A. $$%

It follows %
\begin{equation}\label{eq4.15}
A = \bigcap _{n\geq 1}\cl \left(\bigcup_{m\geq n}A_m\right). %
 \end{equation}  %

Let's prove now that %
\begin{equation}\label{eq4.16}
A\subset  \bigcap_{\epsilon > 0}\bigcup_{n\geq 1}\bigcap_{m\geq n}(A_m)_\epsilon .
\end{equation}
For $0<\epsilon < 1/2$ choose $n_0\in \mathbb N$ such that %
$$%
\forall m\geq n_0 \; F_{AA_m}(\epsilon) > 1- \epsilon. $$%
By \eqref{eq4.12} we have %
$$%
\forall m\geq n_0 \; \forall p\in A \; \exists q\in A_m \; F_{pq}(\epsilon) > 1-\epsilon
, $$%
implying %
$$%
A \subset \bigcap_{m\geq n_0}(A_m)_\epsilon \subset %
 \bigcup_{n\geq 1}\bigcap_{m\geq n}(A_m)_\epsilon
$$%
Again, since $0< \epsilon < 1/2$ was arbitrarily chosen, we get \eqref{eq4.16}.

Finally, prove that %
\begin{equation}\label{eq4.17}
B := \bigcap_{\epsilon > 0}\bigcup_{n\geq 1}\bigcap_{m\geq n}(A_m)_\epsilon \subset %
\bigcap _{n\geq 1}\cl \left(\bigcup_{m\geq n}A_m\right) =: C.
\end{equation}
If   $p\in B$  then %
$$ %
\forall \epsilon , \; 0 <\epsilon < 1, \; \exists n_0(\epsilon )\; \forall m\geq
n_0(\epsilon ) \; p\in (A_m)_\epsilon . $$ %
For $n\geq 1$ letting $m = \max \{ n, n_0(\epsilon )\}$ we have %
$$%
p\in (A_m)_\epsilon \subset  \left(\bigcup_{m'\geq n}A_{m'}\right)_\epsilon .$$%

 We have obtained %
 $$%
 \forall n\geq 1\; \forall \epsilon > 0 \; p\in \left(\bigcup_{m'\geq n}A_{m'}\right)_\epsilon,
 $$%
 implying %
 $$%
 \forall n\geq 1 \quad  p\in \cl\left(\bigcup_{m'\geq n}A_{m'}\right), $$%
 so that %
 $$%
 p\in \bigcap_{n\geq 1}\cl\left(\bigcup_{m'\geq n}A_{m'}\right).$$

Combining now \eqref{eq4.15}, \eqref{eq4.16} and \eqref{eq4.17} we obtain
\eqref{eq4.11}.
\end{proof}

Now we shall prove that the family $P_{fc}(L)$ of all nonempty closed convex subsets of
a complete \v Serstnev random normed space $L$ is
 complete  with respect to the probabilistic Pompeiu-Hausdorff
metric $H$.

A {\it \v Serstnev random normed space} (RN space) is a triple
$(L,\nu,\tau)$ where $L$ is a real linear space, $\tau $ is a
continuous triangle function such that $\tau(D^+\times D^+)\subset
D^+,$ and $\nu$ is a mapping $\nu : L\to D^+$ satisfying the
following conditions:

(RN1)\quad $\nu (p) = \epsilon _0 \; \iff \; p = \theta;$

(RN2) \quad $\nu (ap)(x) = \nu (p)(\frac{x}{|a|}) $ \; for \;
$x\geq 0$\; and\; $a\neq 0;$

(RN3)\quad $\nu(p+q) \geq \tau (\nu(p),\nu(q)),  \; p,q \in L.$

If $(L,\nu,\tau)$ is a \v Serstnev RN space  then %
\begin{equation}\label{eq4.18}
\rho (p,q) = \nu (p-q), \; p,q \in L,
\end{equation}
is a random metric on $L$. The topology of $L$ is the strong topology corresponding to
the random metric \eqref{eq4.18}, and $L$ is a metrizable topological vector space with
respect to this topology.  Random normed spaces were defined and studied by A. N. \v
Serstnev \cite{sers62,sers63e,sers64b} (see also \cite[Ch. 15, Sect. 1]{SS}).

The following result holds: %
\begin{theo}\label{th4.9}
Let $(L,\nu,\tau)$ be a \v Sestnev random normed space with
sup-continuous  triangle function satisfying the condition %
\begin{equation}\label{eq4.19}
\tau (F,G)(x) \geq \sup_{t\in [0,1]} \min \{F(tx),G((1-t)x)\}, \;
\end{equation}
for $x \geq 0\;$ and $ \; F,G \in D^+.$

Then the family $P_{fc}(L)$ of all nonempty closed convex subsets
 of $L$ is closed in $P_f(L)$ with respect to the probabilistic
 Pompeiu-Hausdorff metric $H$, hence complete if the random normed space $L$ is
 complete.

  If $L$ is complete then the family $P_{kc}(L)$ of all nonempty
 compact convex subsets of   $L$ is complete with respect to the probabilistic Pompeiu-Hausdorff
 metric.
 \end{theo}
 \begin{proof}
  Observe first that if  the set $A\subset L$ is convex then the
  set $A_\epsilon $ is convex too.

  Indeed, let $q_1,q_2 \in A_\epsilon$ and $ t_1,t_2 > 0,; t_1+t_2
  = 1.$ If $p_1,p_2 \in A$ are such that $\;
  \nu(p_i-q_i)(\epsilon) > 1-\epsilon,\; i=1,2, $ then
  $t_1p_1+t_2p_2 \in A$ and, by \eqref{eq4.19} and (RN2), %

  $$
    \begin{aligned}
  & \nu(t_1p_1+t_2p_2-(t_1q_1+t_2q_2))(\epsilon) \geq \\
  & \geq \min \{\nu(t_1(p_1-q_1))(t_1\epsilon),
  \nu(t_2(p_2-q_2))(t_2\epsilon)\} \\
  & = \min\{\nu(p_1-q_1)(\epsilon),\nu(p_2-q_2)(\epsilon)\} > 1-\epsilon,
\end{aligned}
$$ %
showing that $t_1q_1+t_2q_2\in A_\epsilon.$

Let now $(A_n)$ be a sequence of nonempty closed convex subsets of
$L$ converging to $A\in P_f(L)$ with respect to $H$. By
Proposition \ref{p4.7} %
$$ %
A = \bigcap_{\epsilon > 0}\bigcup_{n\geq 1}\bigcap_{m\geq
n}(A_m)_\epsilon . $$%
Since each $A_m$ is convex, the same is true for $(A_m)_\epsilon$,
as well as for %
$$%
B_{n,\epsilon} = \bigcap_{m\geq n}(A_m)_\epsilon, \; n=1,2, .... %
$$ %
The union of the increasing sequence $B_{1,\epsilon} \subset
B_{2,\epsilon} \subset ... $ of convex sets will be convex too, so
that their intersection for all $\epsilon >0$ is a convex set.

The assertion concerning the family $P_{kc}(L)$ of all nonempty compact convex subsets
of $L$ follows from Theorem \ref{th4.5} and the first assertion of the theorem.
\end{proof}

\def\cprime{$'$}
\providecommand{\bysame}{\leavevmode\hbox to3em{\hrulefill}\thinspace}


\begin{thebibliography}{100}

\bibitem{Aub}
Jean-Pierre Aubin, {\em Mutational and Morphological Analysis,} Birkh\" auser Verlag,
Boston-Basel-Berlin, 1999.

\bibitem{beg96}
Ismat Beg and Rashid Ali, {\em Some properties of the Hausdorff distance in
probabilistic metric spaces,\,}  Demonstratio Math. {\bf 29} (1996), no.~2,
243--249.

\bibitem{Cho}
Shih-Sen Chang and Yeol~Je Cho, {\em Nonlinear Operators in Probabilistic
  Metric Spaces,\/} Nova Science Publishers, Inc, New York, 2001.

\bibitem{IC}
Gheorghe Constantin and Ioana Istr{\u{a}}{\c{t}}escu,  Elements of
  Probabilistic Analysis and Applications, Editura Academiei, Bucharest, 1981,
 (Romanian).

\bibitem{ICo}
\bysame, {\em Elements of Probabilistic Analysis with Applications,\/} Editura
  Academiei, Bucharest, 1989, (Translated from the Romanian by Victor Giurgiu\c
  tiu).

\bibitem{egb68}
Russell~J. Egbert, {\em Products and quotients of probabilistic metric
  spaces,} Pacific J. Math. {\bf 24} (1968) 437--455.

\bibitem{Fr}
M.~Fr\'{e}chet, {\em Recherches th\' eoriques modernes sur le calcul des
  probabilit\' es,\/} Gauthiers-Villars, Paris, 1936.

\bibitem{H}
Olga Had{\v{z}}i{\'c}, {\em Fixed Point Theory in Probabilistic Metric
  Spaces,\/} Serbian Academy of Sciences, Novi Sad, 1995.

\bibitem{HP}
Olga Had{\v{z}}i{\'c} and Endre Pap, {\em Fixed Point Theory in Probabilistic
  Metric Spaces,} Kluwer Academic Publishers, Dordrecht, 2001.

\bibitem{Hand}
Shouchuan Hu and Nikolas~S. Papageorgiou, {\em Handbook of Multivalued
  Analysis, {V}ol. {I}  Theory,\/} Mathematics and its Applications, vol. 419, Kluwer
  Academic Publishers, Dordrecht, 1997.

\bibitem{Is}
V.~Istr{\u{a}}{\c{t}}escu, {\em Introduction to the Theory of Probabilistic
  Metric Spaces and Applications,\/} Editura Academiei, Bucharest, 1974,
  (Romanian).

\bibitem{Kel}
John L.~Kelley, {\em General Topology,} Van Nostrand, New York 1957.

\bibitem{kol-sos98}
J.~Kolumb\' an and A.~So\' os, {\em Invariant sets in Menger spaces,}  Studia Univ.
Babes-Bolyai, Mathematica {\bf 43} (1998), no. 2, 39-48.

\bibitem{kol-sos01}
J.~Kolumb\' an and A.~So\' os, {\em Invariant sets of random variables in complete
metric spaces,}  Studia Univ. Babes-Bolyai, Mathematica {\bf 46} (2001), no. 3, 49--66.

\bibitem{kol-sos02}
J.~Kolumb\' an and A.~So\' os, {\em Selfsimilar random fractal measure using contraction
method in probabilistic metric spaces,} (Preprint), arXiv:math.PR/0202100v1, 2002.

\bibitem{Lo}
M.~Lo{\`e}ve, {\em Probability Theory,\/} {V}ol. {I}, 4th ed., Springer-Verlag, New
  York, 1977.

\bibitem{men42}
K.~Menger, {\em Statistical metrics,}  Proc. Nat. Acad. U.S.A. {\bf 28}
  (1942) 535--537.

\bibitem{Rad}
Viorel Radu, {\em Lectures on Probabilistic Analysis,\/} West Univ., Faculty of
  Mathematics, Seminar on Probabilty Theory and Applications, Timi\c soara,
  1994.

\bibitem{SS}
B.~Schweizer and A.~Sklar, {\em Probabilistic Metric Spaces,\/} North-Holland
  Publishing Co., New York, 1983.

\bibitem{semp92}
Carlo Sempi, {\em Hausdorff distance and the completion of probabilistic
  metric spaces,}  Boll. Un. Mat. Ital. B (7) {\bf 6} (1992), no.~2,
  317--327.

\bibitem{sers62}
A.~N. {\v{S}}erstnev, {\em Random normed spaces. {Q}uestions of completeness,}   Kazan.
Gos. Univ. U\v cen. Zap. {\bf 122} (1962), kn. 4, 3--20
  (Russian).

\bibitem{sers63e}
\bysame, {\em On the concept of a stochastic normed space},  Dokl. Akad. Nauk
  SSSR {\bf 149} (1963) 280--283 (Russian).

\bibitem{sers64}
\bysame, {\em On a probabilistic generalization of metric spaces,}  Kazan. Gos.
  Univ. U\v cen. Zap. {\bf 124} (1964), kn. 2, 3--11 (Russian).

\bibitem{sers64b}
\bysame, {\em Some best approximation problems in random normed spaces,}  Rev.
  Roumaine Math. Pures Appl. {\bf 9} (1964), 771--789 (Russian).

\bibitem{she69}
H.~Sherwood, {\em On $E$-spaces and their relation to  other classes of probabilistic
metric spaces,}  J. London Math. Soc. {\bf 44} (1969), 441-448.

\bibitem{sib}
David~A. Sibley, {\em A metric for weak convergence of distribution
  functions,}  Rocky Mount. J. Math. {\bf 1} (1971) 427--430.

\bibitem{tar76}
Robert~M. Tardiff, {\em Topologies for probabilistic metric spaces,}  Pacific
  J. Math. {\bf 65} (1976)  233--251.

\bibitem{wa43}
A.~Wald, {\em  On statistical generalizations of metric spaces,}  Proc. Nat.
  Acad. U.S.A. {\bf 29} (1943) 196--197.

\end{thebibliography}
\end{document}